\documentclass[11pt]{article}
\usepackage{amssymb}
\usepackage[top=23mm,bottom=23mm,left=23mm,right=23mm]{geometry}
\title{Lifchitz tail and sausage asymptotics   for stable processes in the Poissonian environment on the Sierpi\'{n}ski gasket}
\author{ \begin{tabular}{ccc}Dorota Kowalska &
 \phantom{and }  & Katarzyna Pietruska-Pa\l{}uba
\\
 Warsaw School of Economics && Institute of Mathematics\\
al. Niepodleg\l{}o\'sci 162 && University of Warsaw\\
02-554 Warsaw, Poland && ul. Banacha 2\\
 && 02-097 Warsaw, Poland
 \end{tabular}}

\newcommand{\barint}{
         \rule[.036in]{.12in}{.009in}\kern-.16in
          \displaystyle\int  }

\def\a{{\alpha d_{w}/2}}
\def\n{{\mathbb{N}}}

\usepackage{color}

\begin{document}

\newtheorem{theo}{\bf Theorem}[section]
\newtheorem{coro}{\bf Corollary}[section]
\newtheorem{lem}{\bf Lemma}[section]
\newtheorem{rem}{\bf Remark}[section]
\newtheorem{defi}{\bf Definition}[section]
\newtheorem{ex}{\bf Example}[section]
\newtheorem{fact}{\bf Fact}[section]
\newtheorem{prop}{\bf Proposition}[section]
\newtheorem{oq}{\bf Open question}

\makeatletter \@addtoreset{equation}{section}
\renewcommand{\theequation}{\thesection.\arabic{equation}}
\makeatother

\newcommand{\ckM}{{\cal K}^{\langle M\rangle}}
\newcommand{\ckz}{{\cal K}^{\langle 0 \rangle}}
\newcommand{\ds}{\displaystyle}
\newcommand{\ts}{\textstyle}
\newcommand{\ol}{\overline}
\newcommand{\wt}{\widetilde}
\newcommand{\ck}{{\cal K}}
\newcommand{\ve}{\varepsilon}
\newcommand{\vp}{\varphi}
\newcommand{\pa}{\partia l}
\newcommand{\rp}{\mathbb{R}_+}
\date{}

\maketitle

\begin{abstract}We obtain the Lifschitz tail asymptotics for the
integrated density of states of the subordinate $\alpha-$stable processes on the Sierpi\'{n}ski gasket $\mathcal G,$ evolving among killing Poissonian obstacles. Simultaneously, we derive the large-time asymptotics for the volume of the $\alpha-$stable sausage on the gasket.
\medskip

\noindent
2010 {\bf MS Classification}: Primary 60J75, 60H25; Secondary 47D08, 28A80.

\end{abstract}

\section{Introduction}

The purpose of this paper is to obtain the Lifschitz tail for the integrated density of states of  random Schr\"{o}dinger operators based on fractional laplacians on the Sierpi\'{n}ski gasket perturbed by killing Poissonian obstacles.


The integrated density of states (IDS, for short) comes into play in the analysis of large-volume systems, when the properties of the spectra of the infinite-volume hamiltonians are difficult to capture. In classical
cases -- when the hamiltonian is based on the Laplace operator in $\mathbb R^d$ -- it has been thoroughly examined (see e.g. \cite{Car}, \cite{Stol} for a review). For Poissonian-type interaction, its existence and behaviour near zero have been analysed also in some nonclassical cases (hyperbolic space \cite{Szn2}, Sierpi\'nski gasket \cite{kpp}, general nested fractals \cite{Shima}). All these papers were concerned with diffusion processes (whose generators are local operators).

For nonlocal operators, the results are not as abundant.
In the classical case (i.e. that of  L\'{e}vy processes on $\mathbb R^d$) the existence of the IDS in ergodic random environment and the asymptotics of related functionals was investigated in \cite{Oku}.
Recently, the existence of the IDS for subordinate Brownian motions on the Sierpi\'{n}ski gasket evolving in random Poissonian environment has been established  in
\cite{Kal-Pie}. Now we will examine the behaviour of the IDS near zero for subordinate stable processes on the gasket and we
will show that the Lifschitz tail is present in this case, with
exponents reflecting the specific scaling of stable processes.
We will also determine the asymptotics of the stable sausage
on $\mathcal G$ -- the two are closely related. The methods we use rely on the enlargement of obstacles technique designed by Sznitman \cite{Szn1}, but adapted to the nondiffusive setting.

More precisely, for $\alpha\in(0,2),$ we consider the symmetric $\alpha-$stable process on the unbounded
Sierpi\'nski gasket  $\mathcal G \subseteq \mathbb R^2,$ and an
independent Poisson point process on $\mathcal G$ with intensity $\nu \mu,$
where $\mu$ is the Hausdorff measure on $\mathcal G$ in dimension $d_f=\frac{\log 3}{\log 2}$ and $\nu>0$ is given.
The points of the Poisson process are centers of balls with radius
$a$ which we call obstacles. The stable process is killed after
entering one of the obstacles. Typically, the semigroup corresponding to this process is not be trace-class and so the spectrum of its generator
may not be discrete.  To get hold on some properties of the spectrum one considers the process in large balls $\mathcal G^{\langle M\rangle}$ of diameter $2^{M}.$  The stable process is then killed when it comes to the obstacle set, or when it jumps out of the set $\mathcal G^{\langle M\rangle}.$ We are interested in the spectra of
the generators of this processes, $\mathcal L^ {M}=\mathcal L^{M,\omega}.$ Now the semigroups become trace-class, so these spectra are pure point and consist of eigenvalues without
nontrivial accumulation points. For each $M$ we build the
empirical measure based on these random sequences of eigenvalues and
normalize them by dividing by the volume of the $\mathcal G^{\langle M\rangle}$. These
empirical random  measures, denoted  by ${l} (M, \omega)$,  converge vaguely, when $M
\rightarrow \infty$, to a deterministic measure ${l}$ on $[0,
\infty)$ which is by definition the integrated density of states.   We will
prove that the IDS fulfills the following property: there exist two constants:
$C>0$ and $D>0$ such that
\begin{equation}\label{intro1}
-C\nu\leq\liminf_{\lambda\to 0} \lambda^{d_f/(\alpha d_w)} \log
l([0,\lambda]) \leq \limsup_{\lambda\to 0} \lambda^{d_f/(\alpha d_w)}
\log l([0,\lambda])\leq -D\nu,
\end{equation}
where $d_w=\frac{\log 5}{\log 2}$ is the walk dimension of $\mathcal G.$
It shows that  the decay of $l$ close to zero is exponential -- faster that for the IDS of
the nonrandom  stable hamiltonian, which is only polynomial.
This is the Lifschitz tail asymptotics, first discovered in 1965 by Lifschitz for disordered quantum systems, and subsequently rigorously proven to hold in various other models (see e.g. the reference list of \cite{BBW}).

To get the desired result, we first derive the asymptotics for the Laplace transform of the measure $l$ (denoted by $L$): there are two positive constans $C_1,D_1$ such that
\[-C_1\nu^{\frac{\alpha}{2}d_w/d_\alpha}\leq \liminf_{t\to\infty}
\frac{\log L(t)}{t^{d_f/d_\alpha}}\leq \limsup_{t\to\infty}\frac{\log L(t)}{t^{d_f/d_\alpha}}\leq -D_1 \nu^{\frac{\alpha}{2}d_w/d_\alpha},\]
where $d_\alpha=d_f+\alpha d_w/2.$ To get (\ref{intro1}),  we employ a Tauberian theorem of exponential type from \cite{Fuk}.

Simultaneously, we establish the asymptotics for the stable sausage
on the gasket in large time: we prove that there exist two constants $C_2,D_2>0$ such that for any $x\in\mathcal G$ one has
\begin{eqnarray*}-C_2\nu^{\frac{\alpha}{2}d_w/d_\alpha}& \leq & \liminf_{t\to\infty}
\frac{\log E_x[\exp{(-\nu\mu(X_{[0,t]}^a))}]}{t^{d_f/d_\alpha}}\\
&\leq& \limsup_{t\to\infty}\frac{\log E_x[\exp{(-\nu\mu(X_{[0,t]}^a))}]}{t^{d_f/d_\alpha}}\leq -D_2 \nu^{\frac{\alpha}{2}d_w/d_\alpha}.\end{eqnarray*}

This it the gasket counterpart of the stable sausage asymptotics from \cite{Don-Var} and of the `Wiener sausage' asymptotics on nested fractals from \cite{kpp1}.

\section{Preliminaries}\label{secprel}
{\bf\em Notation.}
Throughout the paper  $A'$ will denote the complement of a set,
 and $A^\rho$ --- the open $\rho-$neighbourhood of a set.
Generic numerical constants whose actual values are irrelevant for our purposes
will be denoted by the lower case letter $c.$ For  important
constants we will use lower case or capital letters with
subscripts. An `admissible number' is any number of the form
$2^n,$ $n\in \mathbb{Z}.$ When $A\subset \mathcal G$ is a measurable (Borel) set and $(X_t)$ is a stochastic process, then
\begin{eqnarray*}
T_A=\inf \{t\geq 0: X_t\in A\} &\mbox{and}& \tau_A=\inf\{t\geq 0: X_t\notin A\}
\end{eqnarray*}
denote respectively  the entrance and the exit time of $A.$

\subsection{The infinite Sierpi\'{n}ski gasket}

The infinite Sierpi\'nski gasket we will be working on is
defined as a blowup of the unit gasket, which in turn is  the unique fixed point of the hyperbolic iterated function system
in $\mathbb R^2,$ consisting of three maps:
\[
\phantom{}\ \ \ \ \ \ \phi_1(x)=\frac{x}{2},\ \ \ \ \   \phi_2(x)=
\frac{x}{2}+\left(\frac{1}{2},0\right) \ \ \ \ \ \
\phi_3(x)=\frac{x}{2}+\left(\frac{1}{2}, \frac{\sqrt 3}{2}\right). \hfill \]
The unit gasket, $\mathcal G^{\langle 0 \rangle},$ is the unique compact subset of
$\mathbb{R}^2$ such that
\[\mathcal G^{\langle 0\rangle}= \phi_1(\mathcal G^{\langle 0 \rangle})\cup \phi_2(\mathcal G^{\langle 0 \rangle})\cup\phi_3(\mathcal G^{\langle 0 \rangle}).\]
Let $V_{(0)}=\{a_1,a_2,a_3\}=\{(0,0), (1,0),
(\frac{1}{2},\frac{\sqrt 3}{2})\}$ be the set of its vertices. Then
we set:
\[\mathcal G^{\langle n\rangle}=2^n \mathcal G^{\langle 0 \rangle} =((\phi_1^{-1}))^n(\mathcal G^{\langle 0 \rangle}),\] and
\[\mathcal G=\bigcup_{n=1}^\infty \mathcal G^{\langle n \rangle}.\]
Then inductively:
\[ V_{(n+1)}=  V_{(n)}\cup \{2^na_1+ V_{(n)}\}
\cup\{2^na_2 + V_{(n)}\}, \]
\[ V^{\langle 0\rangle} =\bigcup_{n=0}^\infty  V_{(n)}.\]
Elements of $V^{\langle 0\rangle}$ are exactly the vertices of all triangles of
size $1$ that build up the infinite gasket.

The gasket is equipped with the usual Euclidean metric
inherited from the plane.
Observe that in this metric one has
$\mathcal G^{\langle M\rangle}= B(0, 2^M).$ The set $\mathcal{G}$ enjoys the scaling property
\[
2{\cal G}={\cal G},
\]

By $\mu$ we denote the Hausdorff measure on $\mathcal G$ in
dimension $d_f=\frac{\log 3}{\log 2},$  normalized to have
$\mu(\mathcal G^{\langle 0 \rangle})=1.$ The number $d_f$ is called the fractal
dimension of $\mathcal G.$
The measure $\mu$ is a $d_f-$measure, i.e.
there exist two positive constants $a_1,$ $a_2$
 such that for $r>0,$ $x\in\mathcal G$
\begin{equation}\label{eqmeasure}
a_1r^{d_f}\leq \mu(B(x,r))\leq a_2 r^{d_f},
\end{equation}
 Another characteristic number of
$\mathcal G$ is its walk dimension, denoted $d_w(\mathcal G)$ of just $d_w.$  We have $d_w=\frac{\log 5}{\log 2}.$   The spectral dimension of
$\mathcal G$ is by definition $d_s=\frac{2d_f}{d_w}.$

\subsection{The Brownian motion}
On the set $\cal{G}$ one defines the Brownian motion (see
\cite{Bar1,Bar-Per}), denoted by $(Z_t, P_x)_{t\geq 0, x\in \mathcal G}.$  It is a symmetric,
strong Markov, Feller process with continuous trajectories, whose
distribution  is invariant under local isometries of ${\cal G}.$
It has a transition density with respect to the Hausdorff measure
$\mu,$ denoted by $g(t,x,y).$ It is continuous in all its
variables, symmetric in $x,y$. The following scaling
property holds true:
\begin{equation}\label{scal-bm}
g(t, 2x, 2y)=\frac{1}{2^{d_f}}\,g(\frac{t}{2^{d_w}}, x, y)= \frac{1}{3}\,g(\frac{t}{5},x,y),\;\;
t>0,\;\; x,y\in {\cal G},
\end{equation}
where $d_w$ is the walk dimension of $\mathcal G.$
This transition density satisfies the following subgaussian estimates:
there exist constants $a_3, a_4, a_5, a_6>0$ such that for $t>0, x,y\in \mathcal G$ one has
\[ \frac{a_3}{t^{d_s/2}}\,{\rm e}^{-a_4\left(\frac{|x-y|}{t^{1/d_w}}\right)^{\frac{d_w}{d_w-1}}}
\leq g(t,x,y) \leq  \frac{a_5}{t^{d_s/2}} \, {\rm
e}^{-a_6\left(\frac{|x-y|}{t^{1/d_w}}\right)^{\frac{d_w}{d_w-1}}}.\]

In fact, the process in \cite{Bar-Per} is defined on a two-sided gasket, but it can be `folded' to yield the process on the one-sided gasket we are working with.

\subsection{Stable processes on the gasket, definition and relevant
properties}

Following  \cite{Bog-Sto-Szt,Kum-Che,Sto}, $\alpha-$stable
processes on $\mathcal G$  are defined via subordination.
 Fix $\alpha\in(0,2).$
Let $S_t$ be the $\alpha/2-$stable subordinator, independent of $Z$: the L\'{e}vy process on
$[0,\infty)$ with Laplace transform $\mathbb{E}({\rm
e}^{-uS_t})={\rm e}^{-tu^{\alpha/2}};$ let $\eta_t(u),$ $t>0,$
$u\geq 0$ be the density of  the distribution of $S_t.$ Then we
set
\[X_t= Z_{S_t},\quad t\geq 0.\]  This process is  called the
{\em symmetric $\alpha-$stable process} on $\cal G$. As $P[S_t=0]=0,$
$X$ has symmetric transition density given by
\begin{equation}\label{stablesub}
p(t,x,y):=\int_0^\infty \eta_t(u)\,g(u,x,y)\,{\rm d}u.
\end{equation}

\medskip

\noindent The transition density defined by (\ref{stablesub})
satisfies \cite[Proposition 3.2]{Bog-Sto-Szt}:
\begin{itemize}
\item[(i1)] $p(t,x,y)$ is jointly continuous in $(0,\infty)\times {\cal G}\times {\cal G},$
\item[(i2)] the semigroup of operators $T_t$ with kernels $p(t,\cdot,\cdot)$
 is both Feller and strong Feller,
 \item[(i3)] $T_t$ is strongly continuous on $C_0({\cal G}).$
 \end{itemize}
Moreover, $p(t,x,y)$ fulfills  the following estimate
(see \cite[Theorem 3.1]{Bog-Sto-Szt} or \cite[Theorem 1.1]{Kum-Che}):

\smallskip

\noindent there exists a positive constant $A_0=A_0({\cal
G},\alpha)$ such that for $t>0, x,y\in {\cal G},$ $x\neq y,$
\begin{equation}\label{dens-estimate}
\frac{1}{A_0}\, \min\left(\frac{t}{|x-y|^{d_\alpha}},
  t^{-d_s/\alpha}\right) \leq  p(t,x,y)\leq A_0\,\min\left(\frac{t}{|x-y|^{d_\alpha}},
  t^{-d_s/\alpha}\right),
  \end{equation}
where $d_s={2d_f/d_w}$ and $d_\alpha=d_f+\alpha d_w/2,$
and also \begin{equation}\label{diag-est} \frac{1}{A_0}
t^{-{d_s/\alpha}}\leq p(t,x,x)\leq A_0
t^{-{d_s/\alpha}}.\end{equation}

\smallskip

 From the scaling of the transition  density of the
Brownian motion (\ref{scal-bm}) and the scaling of the density of
the subordinator (see e.g. \cite[Formula 8]{Bog-Sto-Szt})
\begin{equation}\label{scal-sub}
\eta_t(u)= t^{-2/\alpha}\eta_1( t^{-2/\alpha} u), \;\; t,u>0,
\end{equation}
we derive the following scaling property for the $\alpha-$stable
density:
\begin{equation}\label{scal-stable}
p(t, 2x,2y)=  \frac{1}{2^{d_f}}\, p(\frac{t}{2^{(\alpha
d_w/2)}},x,y),\;\;\; t>0,\; x,y\in{\cal G}.
\end{equation}
 Indeed, one can write:
  \begin{eqnarray*} p(t, 2x,2y)
&=&
\int_0^\infty g(u,2x,2y)\eta_t(u)\,{\rm d}u\\
&=& \frac{1}{2^{d_f}}\int_0^\infty
g(\frac{u}{2^{d_w}},x,y)\eta_t(u)\, {\rm d}u\\
&=& 2^{d_w-d_f}\int_0^\infty g(\tilde u,x,y)\eta_t(2^{d_w}\tilde u)\,{\rm d}\tilde u\\
&=& 2^{d_w-d_f}\int_0^\infty 2^{-d_w}
g(u,x,y)\eta_{t/2^{(\alpha d_w/2)}}(u)\,{\rm d}u\\
&=& \frac{1}{2^{d_f}}\, p(\frac{t}{2^{(\alpha d_w/2)}},x,y)=\frac{1}{3} p(\frac{t}{5^{\alpha/2}},x,y).
\end{eqnarray*}

\noindent We will need the following estimate on  exit time
from balls:

\begin{fact} \emph{\bf \cite[Lemma 4.3]{Bog-Sto-Szt}} \label{dwa-jeden}
For each $k>1$ there exists $A_{1}=A_{1}(k)$ such that for $x\in \mathcal{G}$, $r>0$, $y \in B(x,r/k)$ we have
\begin{equation}\label{bb}
P_{y}[\tau_{B(x,r)}<t]\leq A_{1}tr^{-\alpha d_{w}/2}.\end{equation}
\end{fact}
\noindent Inequality (\ref{bb}) for $x=y$ gives the estimate for
the supremum of the process:\begin{equation} \label{bbb} \mbox{
for any } x \in \mathcal{G},\;\;\; P_{x}[\sup _{0 \leq s \leq t}
|X_{s}-X_0|>r] \leq A_1 tr^{-\alpha d_{w}/2}.
\end{equation}

Let $U \subseteq \mathcal{G}$ be a bounded open set. By $T_{t}^{U}$ we
denote the $L^2-$semigroup generated by the process killed on exiting
$U$:  for
 functions  $f \in
L^{2}(\mathcal{G},\mu)$ one has $T_{t}^{U} f(x)=E_{x}[f(X_{t}); t<\tau_{U}].$

\begin{fact}\emph{\bf   \cite[Proposition 3.2]{Bog-Sto-Szt}}
The semigroup $(T_{t}^{U})_{t\geq 0}$ has both Feller and strong Feller
properties.
\end{fact}

The scaling of the transition density results in the following
scaling for the principal eigenvalue $\lambda(U)$ (i.e. the
smallest eigenvalue of the generator of the process killed outside
$U$):
\begin{equation}\label{scaling-lambda}
\lambda(2U)= \frac{1}{2^{\alpha d_w/2}} \lambda(U)=\frac{1}{5^{\alpha/2}} \lambda(U).
\end{equation}

\noindent The harmonic measure of an open set is defined classically.
\begin{defi}Let $U\subset\mathcal G$ be open and nonempty, let $x\in U.$ The $P_x-$distribution of $X_{\tau_U}$ is called the harmonic measure of $U.$
\end{defi}

If $U$ is nonempty and bounded then the distribution of
$X_{\tau_{U}}$ is absolutely continuous with respect to $\mu$ on
$\mbox{int}\, (U')$ (see \cite[p. 178]{Bog-Sto-Szt}). Its density is called the Poisson kernel and
denoted by $P_{U}(x,y)$.
We have the following estimates for the Poisson kernels of balls.
\begin{fact}\emph{\bf  \cite[Proposition 6.4]{Bog-Sto-Szt}}\label{mh}
Let
There exists a  constant $A^0_2>0$ such that for each $k>1$, $x_{0}
\in \mathcal{G}$, $r>0$ and for
$A_2=c(\frac{k+1}{k-1})^{d_{\alpha}}A_2^0,$ $\tilde
A_2=c(\frac{k-1}{k+1})^{d_{\alpha}}A_2^0,$  we have
\begin{equation}\label{bbbb}
P_{B(x_{0},r)}(x,z)\leq A_2r^{\alpha
d_{w}/2}|x-z|^{-d_{\alpha}},\,\,\,\,\,\,\,\,\,\,\, x\in
B(x_{0},r),\,\,\, z\in B(x_{0}, k r)',\end{equation}
\begin{equation}\label{bbbbb}
P_{B(x_{0},r)}(x,z)\geq \tilde A_2 r^{\alpha
d_{w}/2}|x-z|^{-d_{\alpha}}, \;\;\;x\in B(x_0,r/k),\;
z\in\mbox{\rm Int}\,(B(x_0,r)').
\end{equation}
\end{fact}

\section{The integrated density of states for stable processes on the gasket evolving among killing Poissonian obstacles}

 Let $\nu>0$ (the intensity) and $a>0$ (the radius of the obstacles)
be fixed. Consider the Poisson point process $\cal N$ with
intensity $\nu \mu$ on ${\cal G},$  denote by $(\Omega, {\cal M},
\mathbb{Q})$ the probability space it is defined on. A ball with radius $a$ (`an obstacle') is attached at each of the Poisson points.  One denotes:
${\cal N}(\omega)=\{x_i\}_{i\in \n},$ $\mathcal N_a(\omega)=\bigcup_i
\overline{B}(x_i,a),$  $\mathcal N(A)=\#\{x_i\in \mathcal N(\omega): x_i\in A\}.$ The set $\mathcal N_a(\omega)$ is called the obstacle set, and the set
 ${\cal O}(\omega)={\cal G}\setminus {\cal N}_a(\omega)$ -- the free open set.
 We assume that the stable process and the
Poisson process are independent. The stable process evolves in $\mathcal O(\omega)$ and is killed
when it jumps to the obstacle set ${\cal N}_a(\omega).$

To define the integrated density of states for such a system, one considers the stable process on a bounded gasket $\mathcal G^{\langle M\rangle},$ $M=1,2,...,$  killed when
 it enters the obstacle set, or when is
jumps out of the interior of ${\cal G}^{\langle M \rangle}.$  Formally
speaking, such a process should be denoted by $X^{(M,{\cal
N}_a)},$ but for the sake of notation we will denote it just
by $X.$ It can be realized in the space of c\`adl\`ag functions and
its transition density (with respect to the Hausdorff measure on
$\cal G$) can be expressed by the usual Dynkin-Hunt formula
\[
p^{M,\omega}(t,x,y)= p(t,x,y)- \mathbb{E}_x[p(t-T_{M,\omega},
X_{T_{M,\omega}},y) {\bf 1}_{T_{M,\omega}<t}],\]  where \[T_{M,\omega}= \inf\{t\geq 0: X_t\in \bigcup_i \overline B(x_i,a)\mbox{ or } X_t\notin \mbox{Int}\,{\cal
G}^{\langle M \rangle}\}\] denotes the entrance time into the
obstacle set or into the closure of  $({\cal G}^{\langle M \rangle})'.$ The transition density
$p^{M,\omega}(t,x,y)$ has a more convenient representation:
\begin{equation}\label{pinned}
p^{M,\omega}(t,x,y)=\left\{
\begin{array}{ll} p(t,x,y)P^t_{x,y}[T_{M,\omega}>t] &\mbox{ for } x,y\in \mbox{Int}\,{\cal G}^{\langle M
\rangle}\cap \mathcal{O}(\omega),\\[2mm]
0 & \mbox{ else. } \end{array}\right.
\end{equation}
 $P_{x,y}^t$ are the bridge
measures: conditional distributions of the process subject to the
condition $X_0=x, X_{t}=y.$  The continuity of $p(\cdot,\cdot,\cdot)$ in time and space variables makes these bridges well-defined \cite[Theorem 1]{Cha-Uri}; we also refer to that paper for more information on Markovian bridges.   Similarly as in \cite[Proposition 4.2]{Oku},  we see that the expression (\ref{pinned})  defines
$\mathbb{Q}$--a.s. a transition density which is symmetric in $x,y.$

In virtue of the representation (\ref{pinned}) and the estimate
(\ref{diag-est}), we see that the semigroup on $L^2({\cal G
}^{\langle M \rangle},\mu),$ associated with kernels
$p^{M,\omega}$,  denoted $(T_t^{M,\omega})_{t\geq 0},$ consists of self-adjoint  trace-class operators,
so its generator $\mathcal L^{M,\omega}$ is self-adjoint
and has pure point spectrum consisting of nonnegative eigenvalues
without accumulation points:
\begin{equation}\label{spectrum}
0\leq \lambda_1(M,\omega) \leq \lambda_2(M,\omega) \leq ...\leq
\lambda_n(M,\omega)\leq...
\end{equation}
One considers then the empirical measures with atoms at points of these
spectra, normalized by the volume of the sets ${\cal G}^{\langle M
\rangle}:$
\begin{equation}\label{l-measure}
l(M,\omega)= \frac{1}{\mu({\cal G}^{\langle M
\rangle)}}\sum_{n=1}^\infty \delta_{\{\lambda_n(M,\omega)\}},
\end{equation}
 and we are
interested in the asymptotical behavior of those measures as
$M\to\infty.$

\smallskip

As in the classical case, these measures have a vague, nonrandom limit $l.$ This limiting measure is called
{\em the integrated density of states} for the $\alpha-$stable
process, or the $\alpha-$stable integrated density of states
($\alpha-$IDS, for short). More precisely, in the paper \cite{Kal-Pie} we have proven the following.

\begin{theo}\emph{\bf \cite[Theorem 3.3]{Kal-Pie}}\label{th:IDS-existence}
Almost surely with respect to the measure $\mathbb{Q},$ the
measures $l(M,\omega)$ converge vaguely to a nonrandom measure $l$
on $[0,\infty).$ \end{theo}

\medskip

 The key to the method is the following representation of the Laplace
transform of empirical measures (\ref{l-measure}):

\begin{eqnarray}\label{eq:el-em}
L(M,\omega)(t) &=& \int_0^\infty {\rm e}^{-\lambda t}{\rm
d}l(M,\omega)(t)\nonumber\\[1mm]
 & = &\frac {1}{\mu({\cal G}^{\langle M
\rangle})}\sum_{n=1}^\infty
{\rm e}^{-\lambda_n(M,\omega)t}
= \frac{1}{\mu({\cal G}^{\langle M \rangle})}
\,\mbox{Tr}\,T_t^{M,\omega}\nonumber\\[2mm]& = &\frac{1}{\mu({\cal G}^{\langle M
\rangle})} \int_{{\cal G}^{\langle M \rangle}}
p^{M,\omega}(t,x,x)\,{\rm d}\mu(x)\nonumber\\[2mm]
& =& \frac{1}{\mu({\cal G}^{\langle M \rangle})}\int_{{\cal G}^{\langle M \rangle}} p(t,x,x)
 P_{x,x}^t[T_{M,\omega}>t]{\rm d}\mu(x).
\end{eqnarray}

The statement of Theorem \ref{th:IDS-existence} was achieved by  proving that for any $t>0$ the averaged Laplace transforms converge to $L(t):$
\begin{equation}\label{el-em-lim}
\mathbb{E_Q}[L(M,\omega)(t)] =: L_M(t)
\stackrel{M\to\infty}{\longrightarrow}
L(t)
\end{equation}
and that the measure $l$ with Laplace transform $L$ is  the $\mathbb Q-$almost sure vague limit of the measures $l(M,\omega).$ Let us note that for the Brownian motion on fractals the limit in (\ref{el-em-lim}) was monotone increasing,
 what followed from symmetries of the process (see \cite{kpp1}), therefore to get the convergence of $L_M(t)$ when $M\to\infty$ one just had to find an upper bound. For stable processes, such symmetries are no longer
true, we cannot use monotonicity, and the proof of the convergence got substantially more difficult.

\section{Asymptotics for the IDS and the stable sausage}
As indicated in the Introduction,
the behaviour of the IDS when $\lambda\to 0^+$ and the asymptotics of its Laplace transform at $+\infty$  are linked via a Tauberian-type theorem, therefore it is enough to get bounds on $L(t)$ when $t\to\infty.$  In the fractal setting, such bounds were previously obtained for the Brownian motion
on the Sierpi\'{n}ski gasket with Poissonian obstacles \cite{kpp}.
 We also refer to \cite{Shima}  for a direct proof of the Lifschitz tail
 for the Brownian motion on nested fractals with potential interaction. Let us note that the methods we use are also suitable for determining  the `sausage asymptotics' when $t\to\infty.$ This topic for the Brownian motion on the Sierpi\'{n}ski gasket was previously addressed in \cite{kpp}, and for general nested fractals -- in
\cite{kpp1}.
Similarly to the Brownian motion case,  the $\alpha-$IDS
asymptotics and the $\alpha-$stable sausage asymptotics are
the same up to a constant, although neither seems to be a direct consequence of the other.

The lower and the upper bounds for the Laplace transform  are obtained separately.
The lower bound estimate is easier: adjusting the ideas
from \cite{kpp1} to the case where one cannot use the monotonicity
of  expressions approximating $L(t),$ we get the desired result by
imposing some additional conditions on the process and on the cloud.
The proof
of the matching upper bound  uses a non-diffusive counterpart of
Sznitman's theorem \cite[Theorem 1.3]{Szn1}, obtained in \cite[Theorem 1]{Kow1}
(see also Theorem \ref{trans} below).

\subsection{The lower bound}

As usual \cite{kpp,Szn1}, to get the lower bound observe that  the event $\{T_{M,\omega}>t\}$ holds true if the
process stays in a sufficiently large ball up to time $t,$ and no Poisson points are present in the vicinity of this ball. For the sausage estimate, one just picks a large ball centered at the origin. For the IDS estimate, the ball will depend on the starting point.  The estimates are then obtained via semigroup methods.

\begin{theo}\label{dolne}
There exist  constants $C_1,C_2>0$ such that for the
Laplace transform of the $\alpha-$IDS one has:
\begin{equation}\label{dolne-IDS}
\liminf_{t\to\infty} \frac{\log L(t)}{t^{d_f/d_\alpha}} \geq -C_1
\nu^{\frac{\alpha}{2}\,{d_w}/{d_\alpha}}
\end{equation}
and for the $\alpha-$stable sausage  volume one has, for any $x\in\mathcal G,$
\begin{equation}\label{dolne-sausage}
\liminf_{t\to\infty} \frac{\log
E_x[\exp(-\nu\mu(X_{[0,t]}^a))]}{t^{d_f/d_\alpha}} \geq
-C_2 \nu^{\frac{\alpha}{2}\,{d_w}/{d_\alpha}}.
\end{equation}
\end{theo}

\noindent {\bf Proof of (\ref{dolne-IDS}).} Let $t>0$ be fixed. As
 $L(t)=\lim_{M\to\infty}L_M(t),$  it is enough to find a lower bound on $L_M(t),$ independent of $M>M_0(t).$
 Recall that $d_\alpha= d_f+\alpha d_w/2$ and let
 \begin{equation}\label{eq:nn}
M_0=M_0(t)=\left[ \frac{1}{d_\alpha}\frac{\log(t/\nu)}{\log
2}\right]\end{equation}
 ($[x]$ denotes the biggest integer not
exceeding $x$). This is the unique integer for which
\begin{equation}\label{eq:nn1}
2^{M_0}\leq \left(\frac{t}{\nu}\right)^{1/d_\alpha} <2^{M_0+1}.\end{equation}

Assuming $M>M_0,$ write
\begin{eqnarray*}
L_M(t)
&=& \frac{1}{\mu({\cal G}^{\langle
M\rangle})}\int_{{\cal G}^{\langle M\rangle}} p(t,x,x)
P^t_{x,x}\otimes\mathbb{Q}[T_{M,\omega}>t]\,{\rm d}\mu(x)\\[2mm]
&= & \frac{1}{\mu({\cal G}^{\langle M\rangle})} \sum_{\mathcal T}\int_{\mathcal T}
p(t,x,x)
P^t_{x,x}\otimes\mathbb{Q}[T_{M,\omega}>t]\,{\rm d}\mu(x),\nonumber
\end{eqnarray*}
where the sum is taken over  all the triangles of size $2^{M_0}$
building $\mathcal G^{\langle M\rangle}.$ These triangles have disjoint interiors
and there are $3^{M-M_0}$ of them. Choose $\mathcal T$ to be one
of those triangles. When $x\in \mathcal T,$ then the event
$\{T_{M,\omega}>t\}$ will hold when the process stays in $\mathcal
T$ up to time $t$ and there are no obstacles in $\mathcal T^a.$
Consequently, we have:
\begin{eqnarray}\label{eq:nowe1}
L_M(t) &\geq & \frac{1}{\mu(\mathcal G^{\langle M
\rangle})}\sum_{\mathcal T} \int_{\mathcal T
}p(t,x,x)P_{x,x}^t\otimes\mathbb Q [\tau_{\mathcal T}
>t, \mathcal N(\mathcal T^a)=0]\,{\rm
d}\mu(x)\nonumber \\[2mm]
&=&\frac{1}{\mu(\mathcal G^{\langle M \rangle})}\sum_{\mathcal T}
\left[\int_{\mathcal T }p(t,x,x)P_{x,x}^t [\tau_{\mathcal T}
>t]\,{\rm d}\mu(t)\right]\cdot\mathbb Q[\mathcal N(\mathcal T^a)=0]\nonumber \\[2mm]
&=&\frac{1}{\mu(\mathcal G^{\langle M \rangle})}\sum_{\mathcal T}
\mbox{Tr}\, T_t^\mathcal T \cdot {\rm e}^{-\nu\mu(\mathcal T^a)},
\end{eqnarray}
where $(T_t^\mathcal T)_{t\geq 0}$ is the Dirichlet stable semigroup on $\mathcal
T .$ Clearly, $\mbox{Tr}\,T_t^\mathcal T \geq {\rm
e}^{-t\lambda(\mathcal T)},$ $\lambda(\mathcal T)$ being the
principal eigenvalue of $\mathcal T$ (relative to the stable
process). It is a classical fact (see e.g.
\cite[Theorem 3.4]{Che-Son}) that
\[\lambda(\mathcal T) \leq (\lambda^{BM}(\mathcal T))^{\alpha/2}\]
where $\lambda^{BM}(\mathcal T)$ is the principal Brownian
Dirichlet eigenvalue of $\mathcal T.$ From symmetry properties of
the Brownian motion on the gasket we see that for any triangles
$\mathcal T, \mathcal T'$ appearing in the sum above one has
$\lambda^{BM}(\mathcal T)=\lambda^{BM}(\mathcal T')=\lambda^{BM}
(\mathcal G^{\langle M_0\rangle}).
 $ We also  have $\mu(\mathcal T^a)\leq
\mu(\mathcal T)+ ca^{d_f}, $ with $c$ -- a numerical constant.
Inserting these bounds into (\ref{eq:nowe1}) we get
\begin{eqnarray*}
L_M(t) &\geq & \frac{1}{\mu(\mathcal G^{\langle M\rangle})} \,
3^{M-M_0}{\rm e}^{-t(\lambda^{BM}(\mathcal G^{\langle
M_0\rangle}))^{\alpha/2}}\cdot {\rm e}^{-\nu(3^{M_0}+ca^{d_f})}
\\[2mm]
&=& \frac{1}{3^{M_0}}{\rm e}^{-t(\lambda^{BM}(\mathcal G^{\langle
M_0\rangle}))^{\alpha/2}-\nu(3^{M_0}+ca^{d_f})}.
\end{eqnarray*}

The scaling of the Brownian principal eigenvalue
($\lambda^{BM}(2U)=\frac{1}{5}\lambda^{BM}(U)$) gives
$\lambda^{BM}(\mathcal G^{\langle M_0\rangle})= \frac{1}{5^{M_0}}
\lambda^{BM}(\mathcal G^{\langle 0 \rangle}). $ Moreover, as
$3^{M_0}=2^{d_fM_0}$ and $5^{M_0}=2^{d_wM_0},$ from (\ref{eq:nn1})
we obtain
\begin{eqnarray*}
5^{M_0}>\displaystyle\frac{1}{5}\,\left(\frac{t}{\nu}\right)^{d_w/d_\alpha} &\mbox{ and } &
3^{M_0}\leq\left(\frac{t}{\nu}\right)^{d_f/d_\alpha},
\end{eqnarray*}
so that
\begin{eqnarray*}
t(\lambda^{BM}(\mathcal G^{\langle M_0 \rangle}))^{\alpha/2}+\nu 3^{M_0}&\leq&
t^{d_f/d_\alpha}\nu^{\frac{\alpha}{2} d_w/d_\alpha}\left({5^{\alpha/ 2}}(\lambda^{BM}(\mathcal G^{\langle 0 \rangle}))^{\alpha/2}+1\right)\\
&=:& C_1 t^{d_f/d_\alpha}\nu^{\frac{\alpha}{2} d_w/d_\alpha},
\end{eqnarray*}
i.e.
\begin{equation}\label{eq:aa}L_M(t)\geq {\rm e}^{- C_1t^{d_f/d_\alpha}\nu^{\frac{\alpha}{2} d_w/d_\alpha}}\cdot{\rm e}^{-\nu ca^{d_f}}\cdot\left(\frac{\nu}{t}\right)^{d_f/d_\alpha}.
\end{equation}
This bound is independent of $t$ and valid as long as $M>M_0(t),$ therefore it is also true for $L(t)=\lim_{M\to\infty} L_M(t).$  Taking
the logarithm, dividing by $t^{d_f/d_\alpha}$ and passing to the limit we obtain (\ref{dolne-IDS}).

\medskip

\noindent {\bf Proof of (\ref{dolne-sausage}).} Let  $t>0$ be given and large enough to have
 $x\in\mbox{Int}\, \mathcal G^{\langle M_0\rangle},$ where $M_0=M_0(t)$
  is introduced in (\ref{eq:nn}).
 Recalling that we have denoted $\mathcal O(\omega)= \mathcal G\setminus \mathcal N_a(\omega),$  we can write
\[E_x\left[\exp[-\nu\mu(X_{[0,t]}^a)]\right]=  P_x \otimes \mathbb{Q}
 [\tau_{\mathcal O(\omega)}>t].\]
 As in the proof of (\ref{dolne-IDS}) we observe that  the event $[\tau_{\mathcal O(\omega)}>t]$ will hold if the process stays in the  ball
 $B(0,2^{M_0})=\mathcal G^{\langle M_0 \rangle}$ up to time $t$ and the $a-$vicinity of this ball receives no Poisson points.    It follows
\[E_x\left[\exp[-\nu\mu(X_{[0,t]}^a)]\right]\geq  \exp[-\nu\mu(({\cal G}^{\langle M_0\rangle})^a)] P_x[\tau_{{\cal
G }^{\langle M_0\rangle}}>t].\] Now write $x=2^{M_0}y$, then scale down using (\ref{scal-stable}) and get
\[P_x[\tau_{{\cal
G }^{\langle M_0\rangle}}>t]= P_y[\tau_{{\cal G}^{\langle 0
\rangle}}>\frac{t}{5^{M_0\alpha/2}}]\geq  P_y[\tau_{{\cal G}^{\langle 0
\rangle}}>5^{\alpha/2}t^{d_f/d_\alpha}\nu^{\frac{\alpha}{2}\,d_w/d_\alpha}].\]
The rest of the proof goes identically as that of \cite[Theorem 2.1]{kpp1}. \hfill $\Box$

\subsection{ The upper bound}
We intend to use the
Sznitman's
`enlargement of obstacles' method in its
non-diffusive  version from  \cite{Kow1}. The method   works for
processes with compact state-space, so the first ingredient needed in the
proof is the reduction of the problem to a one with a compact
state-space.
Once it is done,
the method relies on replacing the microscopic Poisson obstacles
with bigger balls of `intermediate' size and on controlling the
possible increase of principal eigenvalue when the process is
killed on entering those bigger obstacles.

\medskip

To make the paper self-contained, we briefly describe
the method (Section \ref{sec:method}), then we carry out the reduction to the
compact problem (Section \ref{sec:folding}) and prove the necessary estimates
for the `stable process on the compact set' (Section \ref{sec:properties}).
Finally we enlarge the obstacles -- the conclusion of the proof is much alike
that in the Brownian motion case (Section \ref{sec:upper-fin}).

\subsubsection{Description of the method}\label{sec:method}
The method we are going to use works in the following situation:
\begin{itemize}
\item $(\Xi,  d,\mu)$ is a compact metric measure space equipped with
a doubling probability measure $\mu$ charging all open balls. More
precisely, we assume that  there exist constants $\kappa>0$ and
$R_0>0$ such that for any ball $B(x,r),$ $0<r<R_0$ one has
\begin{equation}\label{doubling}
\mu (\overline{B}(x,r/3))\geq \kappa^{-1}\mu (\overline{B}(x,r)).
\end{equation}
\item $(\xi_t)$ is a symmetric strong Markov, Feller process on $\Xi$ with c\`adl\`ag
trajectories and transition
density $p(\cdot,\cdot,\cdot),$  regular
enough to have well-defined symmetric bridge measures
$P^{t}_{x,y}.$
\end{itemize}
Suppose that points  $x_i\in\Xi,$ $i=1,...,n$ are given, together with positive
constants $a> 0,\epsilon>0.$ Closed balls $ \overline{B}(x_i,a\epsilon)$  $i=1,...,n,$  are considered fixed and we call them `obstacles' -- the process is killed when it enters one of
those balls.

\medskip

Points $x_i$ are labeled `good' or `bad' according to
the following rule.

\begin{description}
\item Let $R>0,$  $b\gg a$ and $\delta>0$ be given.
A point $x_{i}$  is called \underline{$(R,b,\delta)-${\em good}}, or
just {\em good}, if for every set $F=B(x_{i},10\epsilon bR^{l})$
such that $10\epsilon b R^{l} \leq R_{0},\,\, l \in \mathbb Z_+ ,$ we have
$$\mu (\bigcup
_{j=1}^{N}B(x_{j},b\epsilon)\cap F)\geq \frac{\delta}{\kappa}\mu
(F).$$
\item Otherwise, $x_i$ is called \underline{$(R,b,\delta)-${\em bad}} (or just {\em bad}).
\end{description}

Denote  $\Theta:=\Xi\setminus  \bigcup_{i=1}^n \overline{ B}(x_i,a\epsilon)$ and for
$b>a$,  $\Theta_b:=\Xi\setminus \bigcup_{i:x_i-\mbox{\scriptsize
good }} \overline B(x_i,b\epsilon),$  then by $\lambda_\Theta$
(resp. $\lambda_b$) -- the principal eigenvalue of the generator
of the processes killed on exiting $\Theta$ (resp. on exiting
$\Theta_b$).

Below we list the  properties of the process which need to be established in
order to make the method work. The numbers $a,b,\epsilon,\delta>0$   are fixed (same as above).

\

We assume that there exist an exponent $s>0$ and numbers $R>3,$ $c_1,c_2,c_3,c_4>0$  such that:

\begin{description}
\item[(A1)]   for
all $x,y \in \Xi$  with  $d(x,y)\leq \beta$, where $\beta$ is an
arbitrary number such that $10b\epsilon \leq \beta <
\frac{R_{0}}{R}$, and for every compact set $E$, which satisfies
$$\frac{\mu(E\cap
\overline{B}(y,\beta))}{\mu(\overline{B}(y,\beta))}\geq
\frac{\delta}{\kappa},$$ we have:
\begin{equation}\label{R1} P_{x}[T_{E}<\tau_{B(y,R\beta)}]\geq c_{1};
\end{equation}
\item[(A2)] with $\beta $ -- as above,  whenever  $x \in \Xi$ satisfies
$d(x_{i},x)\leq \epsilon b <R_{0}$ for some $ i \in \{1,...,n\}$, then
one has:
\begin{equation}\label{R22}
 P_{x}(T_{\overline{B}(x_i,a\epsilon)}\leq
{{\epsilon^{s}}\over{2}})\geq  2c_{2};
\end{equation}
\item[(A3)] for all $x,y \in \Xi,$
when $d(x,y)=b\epsilon$ then
\begin{equation}\label{R21}
P_{y}(\tau_{B(x,10R\epsilon b)}<{{\epsilon^{s}}\over{2}})<c_{2},
\end{equation}
where $c_2$ is the constant from {\bf (A2)};
\item[(A4)]
there exists a decreasing function $\phi:(0,
\infty)\rightarrow (0,1]$ such that if $d(x,y)\leq \epsilon r <
R_{0}$ then
 \begin{equation} P_{x}(T_{B(y, \epsilon b)}\leq  \frac{\epsilon^{s}}{2})\geq \phi (r);
\end{equation}
\item[(A5)]
there exists a constant  $c_{3}>0$  such that for all $0<r<R_0$,
$y \in \Xi$ and $x\in B(y,r)$, when $\rho>3r$, then we have
\begin{equation}\label{R4} P_{x}(
\xi_{\tau_{B(y,r)}}\in B(y,\rho)') \leq c_{3}
\left(\frac{r}{\rho}\right)^{s};
\end{equation}
\item[(A6)] for all  $x,y\in \Xi$ one has
\begin{equation} p(1,x,y)\leq
{c_{4}} .
\end{equation}
\end{description}

The main theorem of \cite{Kow1} is the following.

\begin{theo}\label{trans}Let  $a,b,\epsilon$ be as above
 and let
$K,\delta>0$ be given. Suppose that {\bf
(A1)--(A6)} are satisfied with a number $R>3$ that satisfies
\begin{equation}\label{eq:kow-1}
\frac{c_3}{R^s-1}\leq \frac{1}{2} ({\rm e}^{K}
(1+c_4(1+K/\delta)))^{-1}
\end{equation}
 (the exponent $s$ and the constants
$c_3,c_4$ come from the assumptions above).  Then there exists
$$\epsilon_0=\epsilon_0(b,K,\delta,R_0,R,c_{1},c_{2},c_{3},c_{4},
s, \phi (\cdot))>0$$ such that for $\epsilon<\epsilon_0$ one has
\begin{equation}\label{mainestimate}
\lambda_{b}\wedge K \leq \lambda_{\Theta}\wedge K + \delta.
\end{equation}
\end{theo}

In fact, if we can  find a number $R>3$ for which the assumptions {\bf (A1)--(A5)} hold true, and $R'$ is a number that satisfies (\ref{eq:kow-1}), then the assumptions are fulfilled for $\tilde R=R\vee R'$ as well, without change
in other constants.

\subsubsection{The projected process}\label{sec:folding}
We want to have a process on the finite gasket $\mathcal G^{\langle 0 \rangle},$ locally behaving as
the stable process, and with infinite lifetime. To this goal, we
`project' the unrestricted stable process on $\mathcal G$ onto $\mathcal G^{\langle 0 \rangle},$ using the projection $\pi_0:\mathcal{G}\to\mathcal{G}^{\langle 0\rangle}$ from
\cite[Section 5.2]{kpp}.

Projected processes derived from subordinate Brownian motions
were considered in \cite{Kal-Pie} and they were used there for
proving the existence of the density of states for subordinate processes. Stable processes
on $\cal G$ fall within this category. In particular, the projected
$\alpha-$stable process on $\mathcal G^{\langle 0\rangle}$
 is a strong Markov and Feller process with
continuous, symmetric transition density. See \cite[Section 2.2.3]{Kal-Pie}.

 We recall  briefly its definition.

Following \cite{kpp}, we put labels on vertices from $V^{\langle 0 \rangle}$. We have $V^{\langle 0 \rangle}\subset \mathbb{Z}e_1 +\mathbb{Z}e_2,$ where $e_1=(1,0), $ $e_2= (1/2, \sqrt{3}/2).$ Consider the commutative $3-$group $\mathbb{A}_3$ consisting of even permutations of 3 elements $\{u,v,w\},$ i.e. $\mathbb{A}_3=\{id, p_1, p_2\},$ where
$p_1=(u,v,w), $ $p_2=(u,w,v).$ With every point $x=ne_1+me_2$ we associate the permutation $p_1^n\circ p_2^m\in \mathbb{A}_3.$ In particular, a permutation is assigned to every point $x\in V^{\langle 0 \rangle}$ and the label of $x$ is its value at $u,$ i.e.
$p_1^n\circ p_2^m(u).$

After the vertices have been labeled, we define the projection. Every
nonlattice point  $x\in\mathcal{G}\setminus V^{\langle 0 \rangle}$ belongs to exactly one triangle
of size $1,$  and can be written as \[x=x_u u(x)+x_v v(x) +x_w w(x),\]
where $u(x), v(x), w(x)$ are the points of $\Delta_0(x)$ with respective labels $u,v,w,$ and
numbers $x_u,x_v,x_w\in (0,1)$ satisfy $x_u+x_v+x_w=1.$ For such a point we set
\[\pi_0(x):= x_u\cdot u(0) + x_v\cdot v(0) + x_w \cdot w(0),\]
where we have denoted $u(0)=(0,0), $ $v(0)= (1/2, \sqrt{3}/2),$ $w(0)= (0,1).$
When $x\in V^{\langle 0 \rangle},$ then  $x$ itself has a label and we map it to
the vertex of $\mathcal G^{\langle 0 \rangle}$ with corresponding label.

\noindent Then we define
\begin{equation}\label{proj-def}
Z^{\langle 0\rangle}_t:= \pi_0(X_t)
\end{equation}
and we call this process the {\em projected stable process on $\mathcal G^{\langle 0 \rangle}.$}

In analogy to the reflected Brownian motion from \cite{kpp} whose  transition density
is given by
\[
g^{\langle 0\rangle}(t,x,y)= \left\{ \begin{array}{ll}
\sum _{y' \in \pi _{0} ^{-1} (y) } g(t,x,y) & \textrm{if }  x,y \in \mathcal G^{\langle 0 \rangle}, y \not \in  V^{\langle 0 \rangle}\setminus\{(0,0)\},\\[2mm]
 2 \sum _{y' \in \pi _{0} ^{-1} (y) } g(t,x,y) & \textrm{if }   y \in V^{\langle 0 \rangle}\setminus\{(0,0)\}
\end{array} \right.\]
\noindent the projected stable process (\ref{proj-def}) has transition density $p^{\langle 0\rangle}(t,x,y)$ given by:

\begin{equation}\label{proj}
p^{\langle 0\rangle}(t,x,y)= \left\{ \begin{array}{ll}
\sum _{y' \in \pi _{0} ^{-1} (y) } p(t,x,y) & \textrm{if  } x,y \in \mathcal G^{\langle 0\rangle} , y \not \in V^{\langle 0 \rangle}\setminus\{(0,0)\}\\[2mm]
 2 \sum _{y' \in \pi _{0} ^{-1} (y) } p(t,x,y) & \textrm{if }  y \in V^{\langle 0 \rangle}\setminus\{(0,0)\}.
\end{array} \right.
\end{equation}
Probabilities related to the projected process will be denoted by $(P^{\langle 0\rangle}_x)_{x\in\mathcal G^{\langle 0\rangle}}$.

It is immediate to see that
the projection commutes with subordination, i.e.
\begin{equation}\label{a-1}
p^{\langle 0\rangle}(t,x,y)=\int_{0}^{\infty} g^{\langle 0\rangle}(u,x,y)\eta_{t}(u)\, {\rm d}u\end{equation}
and not hard to prove \cite[Lemma 2.4]{Kal-Pie} that for given $t>0$
the series $\sum _{y' \in \pi _{0} ^{-1} (y) } p(t,x,y')$ is uniformly convergent with respect to $x,y\in{\cal G}^{\langle 0\rangle}.$ The
function $p^{\langle 0 \rangle}$ inherits  symmetry and continuity properties
 of $g^{\langle 0 \rangle},$  established in \cite{kpp}.
The following estimate can be deduced from
\cite[Lemma 2.5]{Kal-Pie}: there exists a constant $A_3>0$ such that for $t>0, x,y\in\mathcal G^{\langle 0 \rangle}$ one has:
\begin{equation}\label{eq:bound-p-zero}
p^{\langle 0 \rangle}(t,x,y)\leq A_3(t^{-{2d_f}/{\alpha d_w}}\vee 1).
\end{equation}

The continuity properties of $p^{\langle 0 \rangle}$ yield the Feller property and then the strong Markov property of the projected process.
Also, these conditions are sufficient for defining bridge measures related to the
projected process (see \cite[Theorems 1,2]{Cha-Uri}). The bridge measure relative to the projected process on $[0,t]$, starting from $x$ and conditioned to arrive at point $y$ at time $t$ will be denoted by
$Q^{t,\langle 0 \rangle}_{x,y}.$

The following proposition permits to relate the bridge of the projected process to the bridge of the free process.

\begin{prop}\label{proj-bridge}
\begin{itemize}
\item[(i)]
Let $x,y \in \mathcal G$ be two points in the same $0-$fiber, i.e. $\pi_{0}(x)=\pi_{0}(y)$. Then the measures $\pi_{0}(P_{x})$ and $\pi_{0}(P_{y})$ on $D([0,t], \mathcal G^{\langle 0 \rangle})$  coincide. Moreover for every $z \in \mathcal G^{\langle 0 \rangle}$ and $x,y$ as above we have:
$$\sum_{z' \in \pi_{0}^{-1}(z)}p(t,x,z')=\sum_{z' \in \pi_{0}^{-1}(z)}p(t,y,z').$$
\item[(ii)]
 For $x,y \in\mathcal G\setminus V^{\langle 0 \rangle}$, the image under $\pi_{0}$ on $D([0,t], \mathcal G^{\langle 0 \rangle})$ of the measure
$$\sum_{y' \in \pi_{0}^{-1}(\pi_{0}(y))}p(t,x,y')P_{x,y'}^{t}[\cdot]$$
is equal to $q^{0}(t,\pi_{0}(x),\pi_{0}(y))Q_{\pi_{0}(x),\pi_{0}(y)}^{t,\langle 0\rangle}[\cdot].$
\\ In particular we have, for all $A \in \mathcal{B}(D[0,t], \mathcal G^{\langle 0 \rangle})$
\begin{equation}\label{most}
q^{0}(t, \pi_{0}(x),\pi_{0}(y))Q_{\pi_{0}(x),\pi_{0}(y)}^{t,\langle 0\rangle}[A]=\sum_{y' \in \pi_{0}^{-1}(\pi_{0}(y))}p(t,x,y')P_{x,y'}^{t}[\pi _{0}^{-1}(A)]
\end{equation}
\end{itemize}
\end{prop}

\noindent{\bf Proof.}
These properties follow from their counterparts  for the Brownian motion  \cite[Theorem 3 and Lemma 8]{kpp} and the subordination formula (\ref{a-1}). See also \cite[Lemma 2.6]{Kal-Pie}.
\hfill$\Box$

\subsubsection{Recurrence properties for $\alpha-$stable processes in
fractals}\label{sec:properties} As we would like  to use Theorem
\ref{trans} for $\Xi = \mathcal G^{\langle 0\rangle},$ $\xi=X^{\langle 0 \rangle},$  we need to establish the relevant recurrence
properties of the projected stable process on
$\mathcal{G}^{\langle 0 \rangle}.$

The $x^{d_f}-$Hausdorff measure on $\mathcal {G}^{\langle 0 \rangle},$ as
well as on $\cal G,$ is  a
$d_f-$measure and as such is doubling. We also know  that the projected $\alpha-$stable
process is strong Markov, Feller, symmetric and regular enough to have well-defined
 bridge measures.

\medskip

\begin{prop}\label{miara}
 Let $X_t^{\langle 0 \rangle}$ be the reflected
$\alpha-$stable process on $\mathcal{G}^{\langle 0\rangle}$
defined by (\ref{proj-def}).
 Let the numbers $b>a>0,$ $\epsilon>0,$ and $\delta>0$  be fixed.  Then {\bf
 (A1)-(A6)} are satisfied, with $s=(\alpha d_w)/2$.
\end{prop}

\noindent\textbf{Proof.} We first check assumptions {\bf
(A1)--(A5)} for the  $\alpha-$stable process on the infinite
fractal. Since for every Borel set  $A\subset \mathcal G^{\langle 0 \rangle}$ and for every
$x\in A$ one has  $P^{\langle 0\rangle}_{x}(\tau
_{A}>t)\geq P_{x}(\tau _{A}>t),$ and for  $y\in \mathcal G^{\langle 0 \rangle}\setminus A$
one has
 $P^{\langle 0\rangle}_{y}(T_{A}\leq t)\geq P_{y}(T_{A}\leq t),$
 conditions {\bf (A1)--(A5)} for the process on the unbounded
 fractal will yield respective properties for the processes on the
 unit fractal.

\smallskip

\noindent {\it Proof of \bf (A1)}. For given $y\in\mathcal G,$ suppose
 $E$ satisfies $\frac{\mu(E\cap
\overline{B}(y,\beta))}{\mu(\overline{B}(y,\beta))}\geq
\frac{\delta}{\kappa}$ and  $|x-y|\leq\beta.$ Is $R$ is large enough (say,  $R>10$) we have
\begin{eqnarray}\label{dzi}
P_{x}[T_{E}<\tau_{B(y,R\beta)}]&\geq& P_{x}[T_{E\cap
B(y,\beta)}<\tau_{B(y,R\beta)}]\nonumber\\[1mm]
&\geq & \inf_{u} P_{x}[X_{\tau_{B(y,2\beta)}} \in B( u, \beta),
  X_{\tau_{B(y,2\beta)}+\tau_{B( u,\beta)} \circ \theta_{\tau_{B(y, 2\beta)}}} \in E\cap B(y,
  \beta)],\nonumber\\
  &&
\end{eqnarray}
where the infimum is taken over $\{u\in \mathcal G : 4\beta
<|y-u|<6\beta\}.$
 From the strong Markov property applied at the stopping time $\tau_{B(y,2\beta)}$ we can
 estimate (\ref{dzi}) from below by
$$\inf_{\{u\in \mathcal G : 4\beta <|y-u|<6\beta\}}\Big(  P_{x}[X_{\tau_{B(y,2\beta)}}
\in B(u, \beta)] \inf _{ z\in B(u, \beta)}P_{ z}[X_{\tau_{B(
u,\beta)}} \in E\cap B(y, \beta)]\Big).$$
 For all $u \in B(y,6
\beta)\setminus B(y,4\beta)$ one has $B(u,\beta)\subset
B(y,3\beta)'.$  Using this fact and the explicit estimate on the Poisson kernel  (\ref{bbbbb}) (with $x_0=y,$  $r=2\beta,$ and $k=2$) we get:
\begin{eqnarray*}
P_{x}[X_{\tau_{B(y,2\beta)}} \in B(u, \beta)]&\geq &\tilde A_2 \int _{B(u,
\beta)} \frac{(2\beta)^{\alpha d_{w}/2}}{|x-z|^{d_{\alpha}}}\,
{\rm d}\mu (z)\\ &\geq&c \beta ^{\a}{\beta ^{-d_{\alpha}}}
\mu(B(u,\beta))\geq ca_1=: c_0.
\end{eqnarray*}
Similarly, using additionally the assumption on the measure of
$E\cap \overline{B}(y, \beta)$ (which is the same as the measure of $E\cap B(y,\beta)$),
\begin{eqnarray*}
 \inf _{ z\in B(u, \beta)} P_{
z}[X_{\tau_{B( u,2\beta)}}
 \in E\cap B(y, \beta)]&\geq& \tilde A_2  \inf _{ z\in B(u, \beta)}
  \int _{E\cap B(y, \beta)} (2\beta)^{\a}\frac{1}{|\zeta-z|^{d_{\alpha}}}\,{\rm d}\mu (\zeta) \\
  &\geq&
c\beta ^{\a} \beta ^{-d_{\alpha}} \mu (E\cap B(y,
\beta))\\
&=& c\beta ^{\a}\beta ^{-d_{\alpha}} \frac{\delta}{\kappa}\beta
^{d_{f}}=\frac{c\delta}{\kappa} =:c_0'.
\end{eqnarray*}
Observe that the constants $c_0$ and $c_0'$ {\em do not} depend on
$\beta.$ Therefore
$$P_{x}[T_{E}<\tau_{B(y,R\beta)}]\geq c_0c_0'=:c_{1}. $$ This
completes the proof of {\bf(A1)}.

\smallskip
{\em Proof of {\bf (A2)}.}
 When $|x-x_i|\leq \epsilon b<R_0,$ then (triangle inequality) for $y\in B(x_i,a\epsilon)$ one has
$|y-x|\leq \epsilon (b+a)\leq 2\epsilon b,$ and one can proceed as follows:
\begin{eqnarray*}
P_{x}(T_{B(x_i,a\epsilon)}\leq \frac{1}{2}\,{\epsilon^{\a}})&\geq&
P_{x}(X_{\frac{1}{2}{\epsilon^{\a}}} \in B(x_i,a\epsilon))\\ &=&
\int_{B(x_i,a\epsilon)}p(\epsilon ^{\a}/2,x,y)\,{\rm d} \mu (y)\\ & \geq &
\mu (B(x_i,a\epsilon)) \inf _{y \in B(x_i,a\epsilon)} \frac{1}{A_0} \min
\Big(\frac{\epsilon^{\a}}{2|x-y|^{d_{\alpha}}},(
\frac{\epsilon^{\a}}{2})^{-d_{s}/\alpha}\Big)\\
&\geq& a_1 (a\epsilon)^{d_{f}}  \frac{1}{A_0}\min
\Big(\frac{\epsilon^{\a}}{2(2b\epsilon)^{d_{\alpha}}},(
\frac{\epsilon^{\a}}{2})^{-d_{s}/\alpha}\Big) \\ &=& c\min
\Big({b}^{-d_\alpha},2^{\frac{d_{s}}{\alpha}}\Big)=:c_2.
\end{eqnarray*}

\medskip

\noindent{\em Proof of} {\bf{(A3)}}. Take $x,y\in\mathcal G$ with
$|x-y|=b\epsilon.$ For $c_2$ chosen as above we just find $R>0$ for which
$$P_{y}(\tau_{B(x,10R\epsilon b)}<\frac{1}{2}\,{{\epsilon^{\alpha
d_{w}/2}}})<c_{2}.$$  This can be done, as according to Fact {\ref{dwa-jeden}},
we have
$$P_{y}(\tau_{B(x,10R\epsilon
b)}<\frac{1}{2}\,{{\epsilon^{\alpha d_{w}/2}}})<(A_1/2)\,
\epsilon^{\alpha d_{w}/2} (10R \epsilon b)^{-\alpha d_{w}/2} =c
(Rb)^{-\alpha d_{w}/2}. $$ Clearly, we can choose $R$ big enough
to make the last quantity  smaller than the previously defined constant
$c_{2}$.

\smallskip

\noindent{\em Proof of }{\bf{(A4)}}. Assume $|x-y|\leq\epsilon
r<R_0,$ so that $B(y,b\epsilon)\subset B(x, (r+b)\epsilon).$  We have the following chain of inequalities:
\begin{eqnarray*}
P_{x}(T_{B(y, b\epsilon)}\leq \frac{1}{2}\,{\epsilon^{\alpha
d_{w}/2}})&\geq& P_{x}(X_{\frac{1}{2}\,{\epsilon^{\alpha
d_{w}/2}}} \in B(y,b \epsilon))\\ &=& \int_{B(y,b \epsilon)}
p(\frac{1}{2}\,{\epsilon^{\alpha d_{w}/2}},x,u){\rm d}\mu (u)\\
&\geq&  a_1 ( b \epsilon)^{d_{f}} \inf _{u \in B(y,b \epsilon)}
p(\frac{1}{2}\,{\epsilon^{\alpha d_{w}/2}},x,u)\\
& \geq & a_1 (b\epsilon)^{d_{f}}  \min\Big(\frac{\epsilon ^{\a}}{
2\epsilon^{d_{\alpha}} (r+b)^{d_{\alpha}}  },
(\frac{\epsilon^{\a}}{2})^{-d_{s}/\alpha}\Big)\\
&=& cb^{d_f}\min ((r+b)^{-d_\alpha}, 2^{d_s/\alpha})=:\phi(r)
.\end{eqnarray*}  The function $\phi$ is nonincreasing (strictly
decreasing for sufficiently big $r$), as required.

\smallskip

\noindent {\em Proof of}  {\bf(A5)}. Assume $\rho>3r,$ $|y-x|\leq r.$ Then the triangle inequality  yields that for $z\notin B(y,\rho)$ one has $|y-z|\leq (4/3) |x-z|.$ From this inequality,  Fact 2.3 (with $k=3$)  and the estimate \[\int_{B(y,\rho)'}|z-y|^{-(d_f+\lambda)} \,{\rm d}\mu(z)\leq c\rho^{-\lambda},\] with $\lambda =\a$
 (see  \cite[Lemma 2.1]{Bog-Sto-Szt})  we obtain
 \begin{eqnarray*}
 P_{x}(
X_{\tau_{B(y,r)}}\in B(y,\rho)')&\leq& A_2 \int_{B(y,\rho)'}
\frac{r^{\a}}{|x-z|^{d_{\alpha}}}\,{\rm d}\mu (z)  \\ &  \leq &   c
r^{\a}\int_{B(y,\rho)'}\frac{1}{|y-z|^{d_{\alpha}}}\;{\rm d} \mu
(z)\\ & \leq & c r^{\a}\rho^{-\a}\leq \frac{c}{3^{\alpha d_w/2}}=:c_3.
\end{eqnarray*}

\medskip

\noindent Finally,  property {\bf{(A6)}} for the projected $\alpha-$stable process follows from (\ref{eq:bound-p-zero}).  \hfill $\Box$

\subsubsection{The upper bound for the Laplace transform and the $\alpha-$stable sausage}\label{sec:upper-fin}
We are ready for the proof of the upper bounds, matching the lower bounds of Theorem
\ref{dolne}.

\begin{theo}\label{gorne}
There exist positive constants $D_1$ and $D_2$ such that for the
Laplace transform of the $\alpha-$IDS one has:
\begin{equation}\label{gorne-IDS}
\limsup_{t\to\infty} \frac{\log L(t)}{t^{d_f/d_\alpha}} \leq -D_1
\nu^{\frac{\alpha}{2}\,{d_w}/{d_\alpha}}
\end{equation}
and for the volume of the $\alpha-$stable sausage  one has: for $x\in \mathcal G,$
\begin{equation}\label{gorne-sausage}
\limsup_{t\to\infty} \frac{\log
E_x[\exp(-\nu\mu(X_{[0,t]}^a))]}{t^{d_f/d_\alpha}} \leq
-D_2 \nu^{\frac{\alpha}{2}\,{d_w}/{d_\alpha}}.
\end{equation}
\end{theo}

\noindent{\bf Proof.} Both (\ref{gorne-IDS}) and
(\ref{gorne-sausage}) are proven similarly as the respective estimates
for the Brownian motion in \cite{kpp}. For clarity, we present the
proof of (\ref{gorne-IDS}) but we skip the other.

Let $t>0$ be fixed. Since for any $t>0$ one has
\(L(t):=\lim_{M\to\infty}L_M(t)=\lim_{M\to\infty}
\mathbb{E_Q}[L(M,\omega)(t)]
,\) it is enough to prove an estimate for $L_M(t)$ which would be
independent of $M.$ As usual, we start with rescaling. Let
$M_0=M_0(t)$ be given by (\ref{eq:nn1}). Writing $x=2^{M_0}y,$ after rescaling we
obtain:
\begin{eqnarray}\label{g-IDS-1}
&&L_M(t)=\frac{1}{\mu(\mathcal{G}^{\langle M-M_0\rangle})}\int_{\mathcal{G}^{\langle M-M_0\rangle}} \frac{1}{2^{M_0d_f}}p(s,y,y)\cdot \nonumber\\[2mm]
&&\hskip 2cm \cdot E^s_{y,y}\left[\exp\left(-\nu 2^{d_fM_0}\mu(X_{[0,s]}^{a/2^{M_0}})\right) {\bf 1}\{\tau_{\mathcal{G}^{\langle M-M_0\rangle}}>s\}\right] {\rm d}\mu(y),\nonumber\\
\end{eqnarray}
where we have denoted $s=t/(2^{M_0\alpha d_w/2}).$

Now we project the process onto $\mathcal G^{\langle 0 \rangle}.$
Starting with the relation
\[p(t,y,y)E_{y,y}^s[\xi]\leq\sum_{y'\in \pi_0^{-1}(\pi_0(y))} p(t,y,y')E_{y.y'}^s[\xi],\]
valid for nonnegative random variables $\xi,$
 splitting the integral over the set $\mathcal{G}^{\langle
M-M_0\rangle}$ into $2^{(M-M_0)d_f}$  integrals over unit cells, then using Proposition
\ref{proj-bridge},  the fact that $\mu(X_{[0,s]}^\rho)\geq \mu
(\pi_0 (X^\rho_{[0,s]}))$ (some volume can be lost in possible
self-intersections of the sausage after the projection), and
neglecting the exit time, we obtain:
\begin{eqnarray}\label{elem-proj}
L_M(t)&\leq & \frac{1}{2^{M_0d_f}}\int_{\mathcal{G}^{\langle 0 \rangle}}
p^{\langle 0 \rangle}(s,y,y)E^{s,\langle 0 \rangle}_{y,y} \left[\exp\left(-\nu 2^{M_0d_f}\mu ((X^{\langle 0\rangle}_{[0,s]})^{a/2^{M_0}}\right)\right]{\rm d}\mu(y)\nonumber\\
\end{eqnarray}
(the bridge measure pertains to the projected process now).

The way $M_0$ was defined gives
\[2^{M_0d_f}\leq \left(\frac{s}{\nu}\right)< 2^{M_0d_f+d_\alpha}\]
so that
\begin{equation}\label{al}
L_M(t)\leq  \frac{2^{d_\alpha}\nu}{s} \int_{\mathcal{G}^{\langle 0
\rangle}} p^{\langle 0\rangle}(s,y,y)E^{s,\langle 0 \rangle}_{y,y}
\left[\exp\left(-2^{-d_\alpha} s \mu((X_{[0,s]}^{\langle 0\rangle})^{
a\nu^{1/(d_f)}/s^{1/(d_f)}})\right)\right]\,{\rm d}\mu(y).
\end{equation}
 The integral
 in (\ref{al}) is equal to the averaged trace  of the semigroup corresponding to the projected
process $X_t^{\langle 0\rangle}$ evolving among (projected and rescaled) killing
obstacles: the intensity of the rescaled Poisson process is  equal
to $\tilde \nu:=2^{-d_\alpha}s$ and the radius of obstacles to $ \tilde
a := a\nu^{1/d_f}/s^{1/d_f}.$
 We can write
 \[(\ref{al}) = 2^{d_\alpha}\nu\,\frac{A(s)}{s},\]
 where $A(s)$ is the averaged trace mentioned:
 \[A(s)= \mathbb{E}^{\widetilde{ \mathbb{Q}}}\int_{\mathcal G^{\langle 0\rangle}}p^{\langle 0 \rangle}(s,y,y)
 P^{s,\langle 0 \rangle}_{y,y}[T_{\mathcal N_{\tilde a}(\omega)}>s]\,{\rm d}\mu(y)\]
($\widetilde{\mathbb{Q}}$ pertains to the rescaled cloud  now).

We now proceed similarly as in the proof of \cite[Lemma 9]{kpp}.
Having proven the recurrence properties {\bf(A1)--(A5)}, we can replace \cite[Theorem 1.4]{Szn1} with
\cite[Theorem 1]{Kow1}, and then obtain \cite[Theorem 1.7]{Szn1}
in the stable case.
In what follows we assume  that the numbers for $\epsilon,b$ are
binary, i.e. of the  form $2^\beta,$ $\beta\in\mathbb{Z}.$ For any
fixed $K,\delta, b>0$ there exists
$\epsilon_0=\epsilon_0(K,\delta,b)$ s.t. for any $\epsilon\leq
\epsilon_0,$ once the radius of obstacles $\tilde a=
a\nu^{1/d_f}/s^{1/d_f}$ is smaller than $a\epsilon_0$ (this
happens when $s$  -- or $t$ -- is large enough) similarly as in
\cite[Formula (77)]{kpp} we get:
\begin{equation}\label{nnum} A(s)\leq c 2^{2(b\epsilon)^{-d_f}}
\exp\Big\{K-\inf_{U\in\mathcal{U}_0}[s(\lambda_0^{\langle 0 \rangle}(U)\wedge
K-\delta)+\tilde\nu (\mu(U)-\delta)]\Big\},\end{equation} where
$\epsilon=(\nu/s)^{1/d_f},$ $\tilde \nu = 2^{-d_\alpha}s,$  $\mathcal{U}_0$
 denotes the collection of all open subsets of $\mathcal G^{\langle 0 \rangle},$ and $\lambda_0^{\langle 0 \rangle}(U)$ is the principal eigenvalue of the reflected
 $\alpha-$stable process on $\mathcal G^{\langle 0 \rangle}$ killed on exiting $U.$ The only difference is that presently we are using part (1) of   \cite[Theorem 1.7]{Szn1}, whereas in \cite{kpp} we were using part (2) of that theorem.

What we get is that, for any $b\gg a,$ $\delta>0, $ $K>0,$ one has
 \begin{equation}\label{nnum1} \limsup_{s\to\infty}
\frac{\log A(s)}{s} \leq \frac{2\ln 2}{b^{d_f}} -\inf_{U\in\mathcal{U}_0}
[(\lambda_0^{\langle 0 \rangle}(U)\wedge K-\delta +2^{-d_\alpha}(\mu(U)-\delta))].
\end{equation} Taking the limits $b\to \infty,$  $\delta\to 0$ and then $K\to\infty$ we
see that
\begin{equation}\label{oops}
\limsup_{s\to\infty}\frac{\log A(s)}{s}\leq -\inf_{U\in\mathcal
{U}_0}[\lambda_0^{\langle 0 \rangle}(U)+2^{-d_\alpha}\mu(U)],
\end{equation}
and as in \cite{kpp}, Lemma 10, we verify that the infimum in
(\ref{oops}) is positive.
Indeed, from (\ref{eq:bound-p-zero}) we get that
that for any $t>0$
\[{\rm e}^{-t\lambda_0^{\langle 0 \rangle}(U)} \leq \int_U p^{\langle 0\rangle, U}(t,x,x)\,{\rm d}\mu(x) \leq A_3 \mu(U) (t^{-2d_f/\alpha d_w}\vee 1),\]
i.e.
\[\lambda_0^{\langle 0 \rangle}(U)\geq -\frac{1}{t}\,\log[A_3\mu(U)(t^{-2d_f/\alpha d_w}\vee 1)].\]
From this estimate it is elementary to see that
we can choose $t$ large enough to guarantee that $\inf_{U\in\mathcal
{U}_0}[\lambda_0^{\langle 0 \rangle}(U)+2^{-d_\alpha}\mu(U)]>0.$

To conclude, observe that the way $s$ was defined gives
\[t^{d_f/d_\alpha}\geq 2^{d_f-d_\alpha} s \nu^{-\frac{\alpha}{2}d_w/d_\alpha}\]
therefore
\[\frac{\log L_M(t)}{t^{d_f/d_\alpha}}\leq 2^{d_\alpha-d_f}\nu^{\frac{\alpha}{2}d_w/d_\alpha}
\frac{\log\left[2^{d_\alpha}\nu\cdot\frac{A(s)}{s}\right]}{s},\] and by passing to the
limit $M\to\infty$ (we can do this as the right-hand side of this
formula {\em does not} depend on $M$) we get the same bound for $L(t).$  Consequently,
\begin{eqnarray*}
\limsup_{t\to\infty}\frac{\log L(t)}{t^{d_f/d_\alpha}}&\leq& \limsup_{s\to\infty}2^{d_\alpha-d_f}\nu^{\frac{\alpha}{2}d_w/d_\alpha}
\frac{\log \left[2^{d_\alpha}\nu\cdot\frac{A(s)}{s}\right]}{s}\\
&\leq&
2^{d_\alpha-d_f}\nu^{\frac{\alpha}{2}d_w/d_\alpha}\inf_{U\in\mathcal
{U}_0}[\lambda_0^{\langle 0 \rangle}(U)+2^{-d_\alpha}\mu(U)].
\end{eqnarray*}
 We

 denote $D_1=2^{d_\alpha-d_f}\inf_{U\in\mathcal
{U}_0}[\lambda_0^{\langle 0 \rangle}(U)+2^{-d_\alpha}\mu(U)] $ and
  (\ref{gorne-IDS}) follows.

\medskip

Inequality (\ref{gorne-sausage}) is proven identically as in
\cite{kpp}, with changes reflecting different scaling, similarly
to the proof of (\ref{gorne-IDS}): after introducing some averaging, the expression estimated can be
compared with the averaged survival time, $B(s),$ of the appropriate
semigroup:
\[B(s)=\mathbb {E}^{\widetilde{\mathbb Q}}\int_{\mathcal G^{\langle 0 \rangle}}P_x[T_{\mathcal N_{\tilde a}(\omega)}>s]\, {\rm d}\mu(x).\]
 Since it is a general fact that $B(s)\leq A(s)$ (see \cite[Formula 1.35]{Szn1}), inequality
(\ref{gorne-sausage}) will follow from the estimates for $A(s)$
proven above. \hfill$\Box$

\subsection{Conclusion. Asymptotics for the $\alpha-$IDS} As in previous articles cited \cite{Szn1,kpp}, Theorems
\ref{dolne} and \ref{gorne} lead to the following estimate,
obtained as an application of the Minlos-Povzner Tauberian Theorem
\cite[Theorem 2.1]{Fuk}.

\begin{theo}\label{final}
There exist two constants: $C=C(D_1)>0$ and $D=D(C_1)>0$ such that
\begin{equation}\label{ff}-C\nu\leq\liminf_{\lambda\to 0} \lambda^{d_s/\alpha} \log
l([0,\lambda]) \leq \limsup_{\lambda\to 0} \lambda^{d_s/\alpha}
\log l([0,\lambda])\leq -D\nu.
\end{equation}

\end{theo}

\noindent
This is the Lifschitz tail asymptotics we intended to prove.


\begin{thebibliography}{99}
\bibitem{Bar1} M. T. Barlow, {\em Diffusion on fractals},
{Lectures on Probability and Statistics, Ecole d'Et\'e de Prob. de
St. Flour XXV --- 1995}, Lecture Notes in Mathematics no. 1690,
Springer-Verlag, Berlin 1998.

\bibitem{Bar-Per}
M.T. Barlow, E.A. Perkins:
\emph{Brownian motion on the Sierpi\'{n}ski Gasket},
Probab. Th. Rel. Fields 79, 543-623, 1988.


\bibitem{BBW} M. Bishop, V. Borovyk, J. Wehr, \emph{Lifschitz Tails for Random Schr\"{o}dinger Operator in Bernoulli Distributed Potentials}, preprint 2014, arXiv:1403.5533.


\bibitem{Bog-Sto-Szt} K. Bogdan, A. St\'{o}s. P. Sztonyk, {\em Harnack inequality
for stable processes on $d-$sets}, Studia Math. {\bf 158} (2) (2003), 163--198.


\bibitem{Car}  R. Carmona, J. Lacroix, {\em Spectral theory of random Schr\"{o}dinger operators}. Probability and its Applications. Birkhäuser Boston, Inc., Boston, MA, 1990.

\bibitem{Cha-Uri} L. Chaumont, G. Uribe Bravo, {\em Markovian bridges: weak continuity
and pathwise construction}, Ann. Prob. {\bf 39}, No. 2 (2011), 609--647.

\bibitem{Kum-Che} Z.Q. Chen, T.  Kumagai, {\em Heat kernel estimates for
stable-like processes on $d$-sets,} Stochastic Process. Appl. 108
(2003), no. 1, 27--62.

\bibitem{Che-Son}
Z.-Q. Chen, R. Song: \emph{Two sided eigenvalue estimates for
subordinate processes in domains}, J. Funct. Anal. 226, 2005,
90--113.

\bibitem{Don-Var} M.D. Donsker, S.R.S. Varadhan, {\em Asymptotics for the Wiener sausage},
 Comm. Pure Appl. Math. 28 (1975), no. 4, 525--565.



\bibitem{Fuk} M. Fukushima, {\em On the spectral distribution of
a disordered system and the range of a random walk}, Osaka J.
Math. {\bf 11} (1974), 73--85.



\bibitem{Kal-Pie} K. Kaleta, K. Pietruska-Pa\l uba,
\emph{Integrated density of states for Poisson-Schr\"{o}dinger perturbations of Markov processes on the Sierpi\'{n}ski gasket}, preprint 2013, arXiv:1310.1027.

\bibitem{Kow1} D. Kowalska, \emph{Lowest Eigenvalue Bounds for Markov Processes with Obstacles},
    Stochastic
Analysis and Applications, {\bf 31} (2013), no. 5, 737-754.

\bibitem{Kus} S. Kusuoka, {\em Dirichlet forms on fractals
and products of random matrices}, { Publ. Res. Inst. Math. Sci} { 25} (1989), 659--680.

\bibitem{Lin} T. Lindstr\o m, {\em Brownian motion on nested
fractals,} Mem. AMS 420 (1990).


\bibitem{Oku} H. Okura, \emph{ On the spectral distributions of certain integro-differential operators with random potential},
 Osaka J. Math. {\bf 16} (1979), 633--666.

\bibitem{kpp} K. Pietruska-Pa{\l}uba,  \emph{The Lifschitz singularity for the density of states on the Sierpinski gasket},  Probab. Th. Rel. Fields {\bf 89} (1991), 1-33.

\bibitem{kpp1}
K. Pietruska-Pa\l uba, {\em The Wiener Sausage Asymptotics on Simple Nested
Fractals}, Stochastic Analysis and Applications, 23:1 (2005), 111-135.

\bibitem{Shima}
T. Shima:
\emph{Lifschitz tails for random Schr\"odinger operators on nested fractals},
Osaka J. Math 29, 1992, 749--770.

\bibitem{Stol} P. Stollmann,  {\em Caught by Disorder. Bound states in random media.} Birkh\"auser, Boston 2001.

\bibitem{Sto} A. St\'{o}s, {\em Symmetric stable processes on $d-$sets},
Bull. Pol. Acad. Sci. Math. {\bf 48} (2000), 237--245.

\bibitem{Szn1} A.S. Sznitman, \emph{ Lifschitz tail and Wiener sausage I}, J. Funct. Anal. {\bf 94} (1990), 223-246.

\bibitem{Szn2} A.S. Sznitman, \emph{ Lifschitz tail and Wiener sausage on hyperbolic space},
Comm. Pure Appl. Math, {\bf 42} (1989), 1033-1065.


\end{thebibliography}
\end{document}